\documentstyle[12pt]{article}
\newcommand{\Bbb}{\mathbf}                       % only for LaTeX
\newcommand{\C}{{\Bbb C}}

\newcommand{\Q}{{\Bbb Q}}

\newcommand{\Z}{{\Bbb Z}}

\newcommand{\cR}{{\cal R}}

\newcommand{\del}{\partial}
\newcommand{\e}{\varepsilon}

\renewcommand{\L}{\Lambda}

\newtheorem{theorem}{Theorem}[section] 
\newtheorem{corol}[theorem]{Corollary}
\newtheorem{cor}[theorem]{Corollary}
\newtheorem{dfn}[theorem]{Definition}

\newtheorem{prop}[theorem]{Proposition}
\newtheorem{rem}[theorem]{Remark}%[section]
\newtheorem{rems}[theorem]{Remarks}%[section]

\newtheorem{eg}[theorem]{Example}%[section]

\newcommand{\be}{\begin{equation}}
\newcommand{\ee}{\end{equation}}
\newcommand{\bea}{\begin{eqnarray}}
\newcommand{\eea}{\end{eqnarray}}
\newcommand{\bean}{\begin{eqnarray*}}
\newcommand{\eean}{\end{eqnarray*}}

\begin{document}
\title{ 
{\bf NONCOMMUTATIVE RATIONAL FUNCTIONS AND 
FARBER'S INVARIANTS OF BOUNDARY LINKS}
}
\author{
Vladimir Retakh,\\ 
{\em Department of Mathematics,} \\
{\em  Rutgers University, \  Piscataway, NJ  08854, USA\/},  \\ 
 \   e-mail: \ {\tt vretakh@math.rutgers.edu} \ \\
 \ \\ 
  Christophe Reutenauer, \\  {\em Departement  de mathematiques, UQAM\/}, \\
 {\em CP 8888 succ Centre-Ville, Montreal, Canada H3C 3P8\/},\\
\ \  \rm e-mail: \ {\tt christo@math.uqam.ca}\\
\ \\and\\
\ \\ Arkady Vaintrob\\
{\em Department of Mathematical Sciences\/,}  \\ 
{\em New Mexico State University, \ Las Cruces, NM 88003-8001, USA\/}, \\
\ e-mail: \ {\tt vaintrob@math.nmsu.edu}
}
\date{}
\maketitle

\bigskip

\hfill {\em To D.~B.~Fuchs on the occasion of his sixtieth birthday}
\bigskip

\begin{abstract}
M. Farber in \cite{Far1} constructed invariants of $\mu$-component  boundary
links with values in the algebra  of noncommutative rational functions.
In this paper we simplify his algebraic
constructions and express them by using noncommutative generalizations
of determinants introduced by Gelfand and Retakh.
In particular, for every finite-dimensional module $N$ over the algebra of
noncommutative polynomials $k\langle x_1,\ldots ,x_{\mu}\rangle $ 
we construct a characteristic rational power series $\chi _N$. 
If $k$ is an algebraically closed field (of arbitrary
characteristic) and $N$ is semi-simple, the series $\chi_N$
determines $N$.
\end{abstract}

\section{Introduction}

Alexander polynomials of knots and links are constructed from the Alexander
module over $\Z[t_1,\ldots,t_\mu]$ (where $\mu$ is the number of components
of the link)  which corresponds to the maximal abelian covering of the link 
complement. They contain essentially all {\em commutative\/} information
about the link (or rather its fundamental group). 

\medskip

Farber in \cite{Far1} constructed invariants of $\mu$-component  boundary
links with values in the algebra  of noncommutative rational
functions. When $\mu=1$ (i.e.\ when link is a knot) this invariant is
equivalent to the Alexander polynomial. 

Farber associated to any $n$-dimensional boundary
link a sequence of rational {\it noncommutative} power
series $\phi _i$, $i=1,\ldots , n$ providing rather
strong link invariants. To compute $\phi _i$ one has
to calculate a finite number of integers (the traces
of certain linear maps acting on the homology of a
Seifert manifold).

The core of Farber's construction is a characteristic rational
power series $\phi _M$ associated to any 
finite-dimensional module $M$ over a certain $k$-algebra
$P_{\mu}$ with $\mu +1$ generators having a subalgebra
isomorphic to the algebra $k\langle x_1,\ldots ,x_{\mu}\rangle $ of
noncommutative polynomials. 
When $k$ is an algebraically closed field of zero characteristic and the
$P_\mu$-module $M$ is semisimple, the series $\phi _M$ determines $M$ up to an
isomorphism. 

In this paper we construct a characteristic rational
power series $\chi _N$ for any finite-dimensional
$k\langle x_1,\ldots ,x_{\mu}\rangle $-module $N$.
When $k$ is an algebraically closed field of {\em arbitrary\/}
characteristic and $N$ is semi-simple, the series $\chi _N$
determines $N$ up to an isomorphism.

If $N$ is also a $P_{\mu}$-module and the actions of
$P_{\mu }$ and $k\langle x_1,\ldots ,x_{\mu}\rangle $
are compatible, the series $\phi _N$ and $\chi _N$
determine each other.

Note that our construction of $\chi _N$ is simpler
than the construction of $\phi _N$ and this gives us 
a hope that this invariant can be applied to a wider class
of links.

Alexander polynomial can be computed
in terms of minors of the Alexander matrix and it would be very 
interesting to have a similar interpretation of its noncommutative
generalization.
  In this paper we  make the first step in the direction of applying
noncommutative determinants to obtain invariants of links.
We show that characteristic series $\chi $ and $\phi $
can be expressed by means of
quasideterminants, a noncommutative generalization of determinants
introduced by Gelfand and Retakh   (see, for example, \cite{GR}).  

\medskip

The first author was suppported in part by
Arkansas Science and Technology Authority and
the second one by NSERC(Canada).

\section{Characteristic functions of finite-dimensional $k\langle X 
\rangle$-modules}

Let $k$ be a field and $X=\{x_1,\ldots,x_\mu\}$ be a set of
noncommuting variables. Denote by $k\langle X\rangle$ (resp.\
$k\langle\langle X\rangle\rangle$) the $k$-algebra of noncommutative
polynomials (resp.\ formal series) in $x_1,\ldots,x_\mu$.
Denote by $X^*$ the free monoid generated by $X$. Then $X^*$ is a 
$k$-basis of $k<X>$. The elements of $X^*$ will be called {\it words}.

\begin{dfn}\rm
The ring $\cR$ of   {\em   noncommutative rational series\/} is defined
as the smallest   $k$-subalgebra of $k\langle\langle X \rangle\rangle$
satisfying the following properties.  

(i) \quad  $k\langle X \rangle \subset \cR$;

(ii) \quad  if $g \in \cR$ and $g$ is invertible in 
$k\langle\langle X \rangle\rangle$ then $g^{-1} \in \cR$.

\end{dfn}

General theory of rational series may be found in \cite {BR}. Here
we remind just basic facts.
\medskip

Rational series are called sometimes {\em noncommutative rational
  functions.\/} They can be characterized in terms of {\em Fox
  derivatives\/}. 

\begin{dfn}\rm
The {\em Fox derivative\/} $\del_i$ on   $k\langle\langle X \rangle\rangle$
is a $k$-endomorphism defined by
$$\del_i(x_jw)=\delta_{ij}w, \quad \del_i 1 =0$$
for any word $w \in X^*$.
\end{dfn}

It satisfies the Leibniz rule: for formal series $f$ and $g$
$$ \partial _i(fg)=(\partial _if)g + \epsilon (f)\partial _ig,$$
where $\epsilon (f)$ means the constant term of $f$.

Fox derivatives turn $k\langle\langle X \rangle\rangle$ into a 
right $k\langle X\rangle$-module. This structure may also be
described by the following construction: the action of polynomial $P$ on
series $S$ is given by $(S\circ P,w)=(S,Pw)$, for
$w\in X^*$, where $(\ ,\ )$ denotes the
canonical pairing between $k\langle \langle X\rangle \rangle $ and
$k\langle
X\rangle$, defined by
$$(S,P)=\sum _{w\in X^*}s_wp_w, \quad \mbox{if}\ S=\sum _{w\in X^*}s_w w,
\quad P=\sum _{w\in X^*}p_w w.$$

In this notations $w\circ x_i=\partial _i w$.
\medskip

\begin{rem} \rm
 Let $\L=k[F_\mu]$ be the group algebra of the free group $F_\mu$ with
$\mu$ generators $g_1,\ldots,g_\mu$. The {\em Magnus embedding\/}
$$\L \hookrightarrow k\langle\langle X \rangle\rangle, 
\quad g_i\mapsto 1+x_i, \quad
g_i^{-1} \mapsto \sum_{n\ge 0} (-x_i)^n
$$
maps $\L$ into $\cR$.
% and the Fox derivatives become the classical free
%derivatives~\cite{Fox}.

By Leibniz rule the image of $\L$ under the Magnus embedding is a
$k\langle X \rangle$-submodule of $\cR$.

The structure of $k\langle X \rangle$-module on $k\langle \langle X \rangle
\rangle$ that we consider here is different from the one obtained by the
canonical embedding $k\langle X \rangle \hookrightarrow k\langle \langle X
\rangle \rangle $.
\end{rem}
\medskip  

Rational series can be characterized in terms of the Fox derivatives. 
The following proposition belongs to 
Sch\"utzenberger~\cite{Sch} and Fliess~\cite{Fl}.

\begin{prop}\label{prop:ratfn}
A series $f$ is rational if and only if the $k$-vector space spanned by all  
Fox derivatives $\del_{i_k}\ldots\del_{i_2}\del_{i_1} f$ of $f$ is
finite-dimensional. 
\end{prop}

Therefore, to each rational function $\chi\in \cR$ there corresponds a
finite-dimentional vector space $M_\chi$ spanned by the Fox derivatives of
$\chi $ of all orders. The space $M_\chi$ has a natural $k\langle X
\rangle $-module
structure given by $x_i \mapsto \del_i, \quad i=1,\ldots,\mu$.
\medskip

Vice versa, to each finite-dimensional $k\langle X \rangle$-module $M$ there
corresponds a rational function $\chi_M$.

\begin{dfn} \rm
Let $M$ be a $k\langle X \rangle$-module, finite dimensional
over $k$, and $u_i \in \mbox{End}_k(M)$ the endomorphism of $M$
corresponding to the action of $x_i \in k\langle X \rangle$.

Define a $k$-algebra homomorphism 
$$\alpha : k\langle X\rangle \to
\mbox{End}_k M, \quad x_i \mapsto u_i.$$

The {\em characteristic function\/} $\chi_M\in k\langle\langle X \rangle
\rangle$ of $M$ is defined by
$$\chi_M=\sum _{w\in X^*}Tr (\alpha (w))w. $$
\end{dfn}

\begin{eg} \rm
Let $k$ be an algebraically closed field and
$\mu =1$, i.e. $M$ is a $k[x]$-module. Let $d=\mbox{dim} (M)$
and $\lambda _i, i=1,\ldots ,d$ be eigenvalues of $\alpha (x)$.
Then
$$\chi _M=\sum _{1\leq i\leq d} \frac {1}{1-\lambda _i x}.$$
\end{eg}
\medskip
 
It is easy to prove the following fact.
\begin{prop}
The characteristic function $\chi_M$ is a rational series 
which is additive for short exact sequences
of  $k\langle X \rangle$-modules.
\end{prop}
\medskip

The following theorem shows that a simple $k\langle X \rangle$-module can be
recovered from its characteristic function.

\begin{theorem} \label{th:1}
Let $k$ be an algebraically closed field and 
$M$ be a simple finite-dimensional $k\langle X \rangle$-module.
Then $M_{\chi_M}$ is isomorphic (as a $k\langle X \rangle$-module) with the
direct sum of $d={\rm dim}(M)$ copies of $M$.
\end{theorem}

We will deduce this result from a theorem of Fliess which
we state now.

Let $M$ be a finite dimensional right $k\langle X \rangle$-module.
Then the dual space $M'$ is a left $k\langle X \rangle$-module.
Let $m\in M$ (resp.\ $\phi \in M'$) generate $M$ (resp.\ $M'$) under the
$k\langle X \rangle$-action.

Define 
$$S=\sum _{w\in X^*}(m(\alpha(w)))\phi w,$$
 where
$\alpha $ is defined as before ($x_i$ acts on the right on
$M'$, $i=1,\ldots,\mu$ and similarly $\phi$ acts on the right on $M$). 

\begin{theorem}[Fliess~\cite{Fl}] Module  $M$ is isomorphic to $M_S$ as
right $k\langle X\rangle$-modules and, under this isomorphism
$m$ (resp.\ $\phi $) corresponds to $S$
(resp.\ to the 
constant term map $\epsilon:
k\langle \langle X \rangle \rangle \to k$).
\end{theorem}

The proof is mechanical: we define the isomorphism by
$m'\rightarrow \sum _{w\in X^*}(m'(\alpha w))\phi w$.
Then a straightforward verification shows that it is well-defined, injective 
and surjective. Note that $k$ needs not to be algebraically closed here.
\bigskip 

\noindent 
{\bf Proof of Theorem~\ref{th:1}}. Let $\{m_1,\ldots,m_d\}$ be
a basis of $M$, and $\{\phi_1,\ldots,\phi_d\}$ be the dual
basis of $M'$.  
    Define $S_i=\sum _{w\in X^*}(m_i(\alpha w))\phi _i$.
Then $\chi_M =\sum _{i=1}^d S_i$.

Let $N$ be a direct sum of $d$ copies of $M$, 
$m=(m_1,\ldots,m_d)\in N$, 
$\phi =(\phi _1,\ldots,\phi_d)\in N'$. Then $\chi_M $
corresponds to the triple $N, m, \phi$ as is
described in the Theorem~\ref{th:1}, so that $M_{\chi_M}$
is isomorphic to $N$, if we can verify the
hypothesis of this theorem.

Since $k$ is algebraically closed, and $M$ is simple
$k\langle X \rangle$-module,  
the subalgebra of
$\mbox{End} _k (M)$ is, by Burnside's theorem, all of
$\mbox{End} _k(M)$. Note that this algebra coincides with 
$\{\alpha P| P\in k\langle X\rangle\}$. Hence, we may find $P\in k\langle 
X\rangle $ such that $m_i(\alpha P)=m_i$ and
$m_j(\alpha P)=0$ for $j\neq i$.
Thus $m(\alpha P)=(0,\ldots,0,m_i,0,\ldots,0)$ which implies
that $m$ generates all of $N$ under this action,
$M$ being simple. For the other hypothesis of Theorem~\ref{th:1}
it is similar: one replaces $M$ by its dual $M'$ on
which the $u_j$'s act simply on the left.

\begin{rem} \rm 
Note that the result is not true if $k$
is not algebraically closed. The smallest example
is $k=\rm R$, $M=ke_1 \oplus ke_2$, and action of $u$ given
by imitating multiplication by complex number $i$,
i.e.
$$e_1u=-e_2,\ \ e_2u=e_1.$$
Then $\chi =\chi_M=2-2x^2+2x^4-2x^6+\ldots$, and
$M_\chi$ is isomorphic to $M$ and not to $M\oplus M$.

There is certainly a version of Theorem~\ref{th:1} when $k$
is not algebraically closed, where the ``arithmetic"
of $k$ plays some role. 
 Also, for effective computational purposes, $k=\Q$ seems to be the best
 field. 
\end{rem}
\medskip

If $M$, instead of being simple, is only semi-simple, 
there is a variant of Theorem~\ref{th:1}.

\begin{theorem}
Let $M$ be a semi-simple module over $k\langle X\rangle$.
If $M=\bigoplus\limits_{r=1}^q M_r,$ where $M_r, \ r=1,\ldots,q$ is simple,
then  $M_{\chi_M}$ is isomorphic to 
$\bigoplus\limits_{r=1}^q\bigoplus\limits_{r=1}^{{\rm dim}(M_r)}M_r.$
\end{theorem}

\medskip

As an application of Theorem~\ref{th:1} we obtain the following important
result.  
\begin{theorem} \label{th:iso} Let $k$ be algebraically closed.
Two finite-dimensional semi-simple 
$k\langle X\rangle$-modules are isomorphic if and only if
their characteristic functions coincide.
\end{theorem}

\bigskip

\section{$P_\mu$-modules}

In his study of boundary links and links modules M.~Farber~\cite{Far1}
introduced an algebra $P_{\mu}$ with $\mu+1$ generators and
defined characteristic functions of finite-dimensional $P_{\mu}$-modules.

Let  $P_\mu$ be a $k$-algebra defined by $\mu+1$
generators  $z,\pi_i, \ i=1,\ldots,\mu,$ and relations 
$$ 
   \pi_i\pi_j=\delta_{ij}\pi_i, \quad \pi_1+\ldots+\pi_\mu =1.
$$

Modules over $P_\mu$ are automatically $k\langle X \rangle$-modules via
the
homomorphism $\delta: k\langle X \rangle \to P_\mu$  defined by $x_i
\mapsto
-z\pi_i$. Denote $\delta (x_i)$ by $\del _i$, $i=1,\ldots ,\mu $.

\begin{prop}
The homomorphism $\delta$ is an embedding.
\end{prop}

\noindent
{\bf Proof.} Indeed, this follows from the two simple facts that, first, the
image of
$\delta$ in $P_\mu$ coincides with the subalgebra $zP_\mu \subset P_\mu$,
and
that, second, this subalgebra is freely generated  by $\mu$ elements $z,
z\pi_1, \ldots, z\pi_{\mu-1}$. Using the identity $\sum_i \pi_i = 1$,
every
word of type  $z^{a_1+1}\pi_{i_1}z^{a_2+1}\pi_{i_2}\ldots z^{a_k}, \quad
a_j \ge 0, \  \ 1\le  i_j \le \mu-1$ can be written as a polynomial in
$\del_i=\delta(x_i)=-z\pi_i$ in a unique way.

\medskip

Recall that $X^*$ is the free monoid generated by 
the alphabet $X=\{x_1,\ldots,x_\mu\}$.
Every word $w\in X^*$ can be uniquely written as $w=x_jw'$ for some $j$ and
$w'\in X^*$. 
Define action of the generators of  $P_\mu$ on $X^*$ as
$$ \pi_i(x_jw) = \delta_{ij}x_j w, \ i,j = 1,\ldots,\mu, \quad z(x_kw)=-w.
$$ 
This gives a $P_\mu$-module structure on the augmentation ideal $k\langle
\langle X \rangle \rangle_+$ of $k\langle \langle X \rangle \rangle $
(i.e.\ the formal series without constant terms).
\medskip

Denote by $\Lambda _0$ the image of the Magnus embedding of
$\Lambda =k[F_{\mu}] \rightarrow k\langle \langle X \rangle \rangle $. The
action of $P_{\mu}$ on $X^*$ defines a $P_{\mu}$-module structure on   
$k\langle\langle X \rangle\rangle /\Lambda _0$. 

%As we have seen in Section 2 the action of  $\del_i=-z\pi_i$  on algebra
%$k\langle\langle X 
%\rangle\rangle$ 
%turns the latter into a right 
%$k\langle X \rangle$-module.   
%Recall that elements $ \del_i = -z\pi_i$ act on $k\langle\langle X
%\rangle\rangle_+$ via the classical Fox derivatives:
%$\del_i(x_jw)=\delta_{ij}w$. We can extend this action to an action on 
%the whole $k\langle\langle X \rangle\rangle$ by setting $\del_i 1 =0$.  
%Note that the operators $\del _i$ generate a  free associative algebra.  

\medskip

\noindent
\begin{prop}
Let $\cR $ be the ring of noncommutative rational series.
Then  ${\cR}/\Lambda _0$ is invariant under $P_\mu$-action on $k\langle
\langle
X \rangle \rangle/\Lambda _0$. 
\end{prop}

\medskip

%\nb{Second example - link modules - mention here ??}
Farber~\cite{Far1} introduced the following
notion of  characteristic function for finite-dimensional
$P_\mu$-modules.
\medskip

\begin{dfn}
Let $A$ be a finite-dimensional $P_\mu$-module. Its characteristic function
is the series
$$ \phi_A = \sum_{k=1}^\mu \sum_\alpha
\mbox{Tr}(\pi_k\del_\alpha)x^\alpha x_k,$$
where $\alpha=(i_1,\ldots,i_p)$, $\del_\alpha= \del_{i_p}\ldots \del_{i_1}$,
$x^\alpha = x_{i_1}\ldots x_{i_p}$. 
\end{dfn}

\begin{prop}
  The characteristic function is a rational series which is additive for
short exact sequences of $P_\mu$-modules.
\end{prop}

\begin{cor} Let $M$ be  finite-dimensional $P_{\mu}$-module
and let $C_i$, be its distinct composition
factors appearing with multiplicities $m_i\geq 1$,
$i=1,\ldots ,n$. Then
$$\phi_M=\sum _{i=1}^n m_i \phi_{C_i}. $$
\end{cor}

\noindent  
\begin{eg} \rm
If $\mu=1$ then 
$$ \phi_A = \sum_{j=1}^n \frac{x}{1+\lambda_j x}, $$
where $n=\mbox{dim}_k A$ and $\lambda_1,\ldots,\lambda_n$ are the eigenvalues
of the operator $z$.
\end{eg}

\medskip

In general, function  $\phi_A$ is an invariant capable of capturing only
{\em semi-simple} information about $P_\mu$-modules. To study non-semisimple 
modules we will need more subtle invariants.
\medskip

\medskip

Farber~\cite{Far1} proved the following result similar to our
Theorem~\ref{th:iso}.
\begin{theorem}\label{th:charp}
Suppose $k$ is an algebraically closed field  of characteristic zero.
Let
$A$ and $B$ be finite-dimensional semi-simple $P_\mu$-modules. Then
$\phi_A=\phi_B$
if and only if $A$ is isomorphic to  $B$.
\end{theorem}
\medskip 

This result follows from the version of our Theorem~\ref{th:1} for
$P_\mu$-modules. 
It shows that the $P_{\mu}$-submodule of ${\cR}/\L _0$
generated by $\chi_A$ is closely related to $A$.

\begin{dfn} \rm
A $P_{\mu}$-module $M$ is called {\em primitive}
if one of the generators  $\pi_1,\ldots , \pi_{\mu}, z$ of $P_\mu$
acts as the identity on $M$ and the other
generators act trivially. 

A $P_{\mu}$-module $A$ is called 
{\em primitive-free} if it has no primitive composition factors.
\end{dfn}

\begin{theorem}[\cite{Far1}]
If $k$ is algebraically closed and $A$ is simple non-primitive $P$-module
of {\rm dim}$A=d$, then
$P_{\mu}\phi_A \subset {\cR}/\L_0$ is $P_{\mu}$-isomorphic to $dA$,
the direct sum of $d$ copies of $A$.
\end{theorem}

\medskip

\begin{rems} \rm
\

(1) \ 
If one takes an alphabet with $\mu +1$ letters,
then the characteristic series $\chi_M$  on this alphabet for a
$P_{\mu}$-module is equivalent to Farber's series; i.e. the knowledge of one 
is equivalent to the knowledge of the other. 

(2) \ It it easy to construct example of simple $P_{\mu}$-modules
which are not $k\langle x_1,\ldots ,x_{\mu}\rangle$-simple
under the homomorphism $x_i\mapsto -z\pi_i$. Let $M$ be a 
three-dimensional  vector space with a basis $e_i$, $i=1,2,3$.
Define an action of algebra $P_3$ on $M$ by setting
$\pi _ie_j=\delta _{ij} e_j, \quad i,j=1,2,3$. 
The action of $z$ in this basis is given by a matrix 
\parbox{1.5in}{$$
\left[ \begin{array}{rrr}
         1&0&1\\
         1&1&2\\
         0&1&1
       \end{array}
\right].
$$}

One can see that $M$ is a simple $P_3$-module under this action, and
the two-dimensional image of $z$ is invariant under action of $\del_i$ for
$i=1,2,3$. 

Therefore, it is not clear whether the main results
about semi-simple $P_\mu$-modules
of this section follow directly from the results obtained in Section~2. 

\end{rems}

\section{Quasideterminants and characteristic functions}

In this section we express characteristic functions constructed
in Sections 2 and 3 via {\it quasideterminants}. 

Quasideterminants were introduced by Gelfand and Retakh (see, for example, 
\cite{GR}) and are defined  as follows. 
Let $A$ be an $m\times m$-matrix over an algebra $R$.
For any $1\le i,j\le m$, let
$r_i(A)$, $c_j(A)$ be the i-th row and the j-th column of $A$.
Let $A^{ij}$ be the submatrix of $A$ obtained by removing
the i-th row and the j-th column from $A$. For a row
vector $r$ let $r^{(j)}$ be $r$ without the j-th entry.
For a column vector $c$ let $c^{(i)}$ be $c$ without the i-th entry.
Assume that $A^{ij}$ is invertible. Then the quasideterminant
$|A|_{ij}\in R$ is defined by the formula
$$
|A|_{ij}=a_{ij}-r_i(A)^{(j)}(A^{ij})^{-1}c_j(A)^{(i)},
$$
where $a_{ij}$ is the $ij$-th entry of $A$.
\medskip
Let $B=(b_{ij})$ be a matrix of order $n$ with formal entries
and $E_n$ a unit $n$-matrix. The following proposition was
proved in \cite{GKLLRT}.

\begin{prop} In the ring $\Z\langle \langle b_{ij} \rangle \rangle$
 of formal series with integer
coefficients generated by noncommuting variables $b_{ij}$ we have
$$ 
|E_n - B|^{-1}_{ii} = 
  1+\sum_{k_1,\ldots,k_p} b_{ik_1}b_{k_1k_2}\ldots b_{k_p i},  
$$
where the sum is over all $1\leq k_1,\ldots ,k_p\leq n$,
$p=1,2,\ldots $.
\end{prop}

\medskip

Let $k$ be a field, $X=\{x_1,\ldots ,x_{\mu}\}$ set of noncommuting
variables, and $M$ be an $n$-dimensional
 $k\langle \langle X\rangle\rangle $-module. 
Choose a basis of $M$ and 
let $A_i$, $i=1,\ldots , \mu,$ be the matrix of the operator $x_i$ 
in this basis.
 
The following proposition gives an expression of 
the characteristic function $\chi_M$ of $M$ via quasideterminants.

\begin{prop} 
The sum 
$$\sum _{i=1}^n |E_n - x_1A_1-\ldots 
-x_{\mu}A_{\mu}|_{ii}^{-1}$$
does not depend on the choice of the basis and is equal to the characteristic
function $\chi _M$ of the $k\langle X \rangle$-module $M$.  
\end{prop}
\medskip

A similar result holds for Farber's characteristic
functions. 

Let  $M$ be a finite-dimensional  $P_{\mu}$-module. Fix a basis of $M$ and
denote by $B_i$, $i=1,\ldots , \mu$  the matrix of the
the operator $(1-z)\pi_i$ in this basis.

\begin{prop}
The sum 
$$\sum _{i=1}^n |E_n - x_1 B_1-\ldots 
-x_{\mu}A_{\mu}|_{ii}^{-1}$$
does not depend on the choice of the basis and is equal to the characteristic
function $\phi_M$ of the $P_\mu$-module $M$.  
\end{prop}

\section{Link modules}

Let $F_\mu$ be a free group with $\mu$ generators $t_i$ and
 $\L=k[F_\mu]$ its group ring. A finitely-generated left $\L$-module $M$ is 
 called {\em link module\/} if $\mbox{Tor}^\L_q(k,M)=0$ for 
all $q\ge 0$, 
 where  $k$ is regarded as a right $\L$-module via the augmentation map.

The following  characterization of link modules allows to consider them as
$P_\mu$-modules.

 \begin{theorem}[Sato~\cite{Sato}]
$M$ is a link module if and only if every element $m$ of $M$ has a unique
representation $$m=\sum_{i=1}^\mu (t_i-1)m_i, \ m_i \in M.$$
 \end{theorem}

 \begin{corol}
In the notations of the theorem, formulas
$\pi_i m=(t_i-1)m_i, \ z m = \sum m_i$, make every link module a
$P_\mu$-module. In addition, $\del_i m = m_i$.
 \end{corol}

Note that every $P_{\mu }$-module is also a 
$k\langle x_1,\ldots ,x_{\mu}\rangle $-module
under the canonical embedding  $\delta: k\langle x_1,\ldots ,x_{\mu}\rangle
\to P_{\mu }$. 
\medskip

\begin{dfn}[\cite{Far2}] \rm
  A finite-dimensional $P_\mu$-submodule $A$ of a link module $M$ is called
  {\em   lattice\/} if $A$ generates $M$ over $\L$.
\end{dfn}

It is easy to prove the following fact.

\begin{theorem}[\cite{Far2}]
  Every link module contains a lattice. Intersection of all lattices in $M$
  is a lattice which is called the minimal lattice of $M$.
\end{theorem}

\begin{dfn}
If $M$ is a link module, we set  $\phi_M = \phi_A$ and
  $\chi_M = \chi_A$, where $A$ is the minimal
lattice in $M$. 
\end{dfn}

\begin{theorem}[Farber~\cite{Far1}]
  Let $M$ and $N$ be semi-simple link modules with $\phi_M=\phi_N$.
If  the field $k$ is algebraically closed and  ${\rm char}(k)=0$, 
then $M$ and $N$ are isomorphic as $P_\mu$-modules.
\end{theorem}

The proof of the theorem is based on Theorem~\ref{th:charp}.

\section{Boundary links}

An $n$-dimensional $\mu$-component link is an oriented smooth submanifold
$\Sigma$ of $S^{n+2}$, where $\Sigma=\Sigma_1\cup \ldots \cup \Sigma_\mu$ is
an ordered disjoint union of $\mu$ submanifolds of $S^{n+2}$, each
diffeomorphic to $S^n$. It is called a {\em boundary link\/} if there exists
an oriented submanifold $V$ of dimension $n+1$ of $S^{n+2}$, such that
$V=V_1\cup \ldots \cup V_\mu$ and $\del V_i=\Sigma_i, \ i=1,\ldots,\mu$. 
If, in addition, each $V_i$ is connected, we say that $V$ is a 
{\em Seifert manifold\/} for $\Sigma$. 
\medskip

The  homology groups of the Seifert manifold $H_*(V;k)$ have several  natural
operations.  The first $\mu$ operations are the projections $\pi_i: H_*(V;k)
\to H_*(V_i;k)$. The last operation $z$ is defined as follows.
For $Y=S^{n+2}\setminus V$ let $I_\pm : V \to Y$ be small shifts in the
direction of positive and negative normals to $V$, respectively. The map
$$ i_{+*} - i_{-*}: H_k(V) \to H_k(Y) $$ 
is an isomorphism for any $k=0,1,2,\ldots$ (see~\cite{Far1}) and we define $z$ 
by $z=(i_{+*} - i_{-*})^{-1}i_{+*}$.

Let $\Sigma$ be a link in $S^{n+2}$ and $X=S^{n+2}\setminus T(\Sigma)$ the
complement of its tubular neighborhood.   

Fix a base point $* \in X$. Then connecting the meridians of $\Sigma$ (small
loops around each component $\Sigma_i$) with $*$ we obtain elements 
$m_1,\ldots,m_\mu \in \pi_1(X,*)$ defined up to conjugation. If $\Sigma=\del
V$ is a boundary link then there is an epimorphism (defined up to
conjugation) $\sigma:  \pi_1(X,*) \to F_\mu$ defined as follows.
If $\alpha$ is a loop in $X$ intersecting $V$ transversally, first in
component $V_{i_1}$,  then  $V_{i_2}$, etc, we set 
$\sigma[\alpha]=t_{i_k}^{\e_k}\ldots t_{i_2}^{\e_2} t_{i_1}^{\e_1}$,
where $\e_i=\pm 1$ is the local intersection index between $\alpha$ and $V$
at $i$-th intersection point.

Consider the covering $\widetilde X \to X$ corresponding to the kernel of
$\sigma$. The group of deck transformations of this covering space is 
$F_\mu$ and therefore the homology $H(\widetilde X; k)$ is a
$\L=k[F_\mu]$-module. For $0\le i \le n$ this module is a link module.

There is a canonical lifting map $S^{n+2}\setminus \Sigma \to \widetilde X$
giving a homomorphism of $P_\mu$-modules $f: H_i(V;k)\to H_i(\widetilde X;k)$.
The image of $f$ is a lattice.  
If $H_i(V;k)$ has no primitive $P_\mu$-submodules, then $f$ is a 
monomorphism and the image of $H_i(V;k)$ is a minimal sublattice of
$H_i(\widetilde X;k)$ (cf.~\cite{Far1, Far2}).  

\bigskip

Applying the constructions of Sections 2 and 3 to the minimal lattices in 
homology groups $H_i(\tilde X, k)$ we obtain sequences of rational
series $\chi _i$ and $\phi _i$, $i=1,\ldots ,\mu$. Thus we have a sequence
of {\it noncommutative invariants} associated to every boundary link. 

These invariants are stronger than the well-known commutative invariants. 
In particular, Farber~\cite{Far1} constructed an example of a link one of
whose Alexander modules  vanishes, but the corresponding characteristic
function is non-trivial.
\medskip

\begin{eg}\rm  For $\mu=1$ the invariants $\phi$, $\chi$,  and the Alexander
polynomial determine each other. Let 
$\Delta _i(t)$ be the Alexander polynomial of $H_i(\tilde X;\Q)$
where $\tilde X$ is an infinite cyclic cover of the complement of
the knot. Then if
$$
\Delta _i(t)=\prod _{j=1}^n (t-\nu _j),\ \ \nu _i\in \C,
$$
the characteristic functions are
$$\chi (x)=\sum _{j=1}^n 1/(1-\lambda_jx)$$
and
$$\phi (x)=\sum _{j=1}^n x/(1+\lambda_jx),$$
where
$$ \lambda_j=1/(1-\nu_j),\ \ j=1,\ldots ,n. $$
\end{eg}

\end{document}